\documentclass{amsart}

\usepackage{geometry}
\usepackage{hyperref}
\usepackage{mathtools}

\theoremstyle{theorem}
\newtheorem*{theorem}{Theorem}

\newcommand{\calM}{\mathcal{M}}
\newcommand{\FF}{\mathbb{F}}
\newcommand{\PP}{\mathbb{P}}
\newcommand{\bfa}{\mathbf{a}}
\newcommand{\bfb}{\mathbf{b}}

\begin{document}

\title{Roots of Polynomials and The Derangement Problem}
\markright{Roots of Polynomials and The Derangement Problem}
\author{Lior Bary-Soroker}
\address{School of Mathematical Sciences, Tel Aviv University, Ramat Aviv, Tel Aviv 6997801, Israel}
\email{
barylior@post.tau.ac.il
}
\author{Ofir Gorodetsky}
\address{School of Mathematical Sciences, Tel Aviv University, Ramat Aviv, Tel Aviv 6997801, Israel}
\email{
ofir.goro@gmail.com
}

\maketitle

\begin{abstract}
We present a new killing-a-fly-with-a-sledgehammer proof of one of the oldest  results in probability which says that the probability that a random permutation on $n$ elements has no fixed points tends to $e^{-1}$ as $n$ tends to infinity. Our proof stems from the connection between permutations and polynomials over finite fields and is based on an independence argument, which is trivial in the polynomial world.
\end{abstract}

\section{Introduction.}
The derangement problem, first studied by Pierre R\'{e}mond de Montmort in 1708, asks what is the probability $P_n$ that a random permutation on $n$ letters has no fixed points. A simple argument using the inclusion-exclusion principle, given by Nicholas Bernoulli in 1713, computes $P_n$ explicitly (namely, $P_n = \sum_{i=0}^{n} \frac{(-1)^i}{i!}$) and, in particular,
\begin{equation} \label{eq:MT}
\lim_{n\to \infty} P_n = e^{-1}.
\end{equation}
If $X_k = \{ \sigma(k) \neq k\}$ is the event that $k$ is not a fixed point, then $P_n = \mathbb{P}(X_{1}\cap \cdots \cap X_{n})$. 
Since $\mathbb{P}(X_k) = 1-\frac{1}{n}$, one gets that $\prod_{k}\mathbb{P}(X_k) = (1-\frac{1}{n})^{n} \to e^{-1}$. 
However, this argument does not imply \eqref{eq:MT} since the events $X_k$ are not independent. 

Our goal is to give a new proof of \eqref{eq:MT} based on an independence argument of polynomials over finite fields. This approach may be considered natural for finite field theorists due to the ancient connection between polynomials, permutations, and integers which goes back at least to Gauss's exact formula; see \eqref{Gauss} below. For further reading on recent progress see the survey papers \cite{BS,KRW,Rudnick}.

We approximate $P_n$ by the probability $P_{n,q}$ that a random uniform monic polynomial of degree $n$ over a finite field $\mathbb{F}_q$ with $q$ elements has no root in $\mathbb{F}_q$ (here $q$ is a prime power): If $n\geq q$, then the events $\tilde X_\alpha = \{ f(\alpha)\neq 0\}$, $\alpha\in \mathbb{F}_q$ are independent  and $\mathbb{P}(\tilde X_\alpha)=1-q^{-1}$. (This is straightforward  from Lagrange interpolation or from the Chinese remainder theorem, see \S\ref{secindep}.)
Thus,
\begin{equation}\label{eq:PT}
P_{n,q} = \mathbb{P}\Big(\bigcap_{\alpha\in \mathbb{F}_q} \tilde X_{\alpha}\Big) = \prod_{\alpha\in \FF_q} \mathbb{P}(\tilde X_\alpha)=(1-q^{-1})^q \to e^{-1}, \quad (n\geq q\to \infty).
\end{equation}
A special case of a theorem of Arratia, Barbour, and Tavar\'{e} \cite[Cor.~5.6]{ABR} says that 
\begin{equation}\label{eq:MB}
|P_n-P_{n,q}|= O\Big(\frac{1}{q}\Big);
\end{equation} 
see \S \ref{proofest} for a simplified  proof of \eqref{eq:MB}. 
Substitute $q=2^{\lfloor \log_2 n \rfloor}$ in \eqref{eq:MB} and take $n\to \infty$ to get by \eqref{eq:PT} that 
\[
\lim_{n\to \infty} P_n = \lim_{n\to \infty} P_{n,q} = e^{-1},
\]
as needed. \qed

\section{Independence argument.}\label{secindep}
Assume $n\geq q$. By Lagrange interpolation, for any subset $S$ of $\FF_q$ and any choice of $(c_\alpha)_{\alpha\in S}$ with $c_\alpha\in \FF_q$, there exists a unique  polynomial $g(X)$ of degree $<\#S$ such that $g(\alpha) = c_\alpha$ for all $\alpha\in S$. 
Thus the monic  polynomials of degree $n$ passing through $(\alpha,c_{\alpha})_{\alpha\in S}$ have the form 
\[
f(X) = g(X) + h(X) \prod_{\alpha\in S} (X-\alpha),
\]
where $h$ is monic and  $\deg h = n-\#S$ and there are $q^{n-\#S}$ of them. (This parametrization also follows from the Chinese remainder theorem.) Thus
\[
\PP\Big(\bigcap_{\alpha\in S} \tilde{X}_\alpha\Big) = \sum_{c_\alpha\neq 0} \PP(f(\alpha)=c_\alpha, \ \alpha\in S) = \sum_{c_\alpha\neq0} \frac{q^{n-\#S}}{q^n} 
= \Big(\frac{q-1}{q}\Big)^{\# S},
\]
which proves independence. Here the sum runs over all tuples of $(c_\alpha)_{\alpha\in S}$ of nonzero elements in $\FF_q$. \qed

\section{Probability measures on the space of partitions.}\label{meas}

Consider the space $\Omega$ of all $n$-tuples $\bfa=(a_1,a_2,\ldots, a_n)$ of nonnegative integers such that $\sum i a_i=n$. We define two probability measures on $\Omega$, one coming from permutations and the other from polynomials over a finite field and compare them.

Let $S_n$ be the symmetric group on $n$ elements. 
Each $\sigma\in S_n$ has a cycle structure which may be regarded as an element $\bfa$ in $\Omega$: We set $a_i$ to be the number of orbits of $\sigma$ of length $i$. Then $\sum_i i a_i = n$. For example, the trivial element corresponds to $(n,0,\ldots,0)$, a transposition to $(n-2,1,0,\ldots,0)$, $(0,\ldots, 0,1)$ an $n$-cycle, etc. The uniform measure on $S_n$ then induces a probability measure $\PP_{S_n}$ on $\Omega$, which by Cauchy's formula is given by   
\begin{equation}\label{eq:ForSn}
\PP_{S_n}(\bfa) = 
\prod_{k=1}^{n} \frac{1}{k^{a_k}a_k!}.
\end{equation}

Let $q$ be a prime power, $\FF_q$ the finite field with $q$ elements, $\FF_q[X]$ the ring of polynomials with coefficients in $\FF_q$, and denote by $\mathcal{M} = \mathcal{M}_{n,q}$ the set of monic polynomials of degree $n$ in $\FF_q[X]$. 
The unique factorization of $f\in \mathcal{M}$ to monic irreducible polynomials, 
\[
f= P_1 \cdots P_k,
\]
defines an element $\bfa$ of $\Omega$; namely, 
\[
a_i = \#\{j : \deg P_j = i\}.
\]
We emphasize that the counting is with multiplicity. Hence $\sum_i i a_i = \deg f= n$. The uniform measure on $\mathcal{M}$ then induces a probability measure on $\Omega$ given by 
\begin{equation}\label{eq:PMnq}
\mathbb{P}_{\mathcal{M}} (\bfa) = 
\frac{\prod_{i=1}^n \pi_q(i,a_i)}{q^n}, 
\end{equation}
where $\pi_q(i,a_i)$ is the number of ways to choose $a_i$ monic irreducible polynomials of degree $i$ (with repetition). If we denote by $\pi_q(i)= \pi_q(i,1)$ the number of monic irreducible polynomials of degree $i$, then $\pi_q(i,a_i) = \binom{\pi_q(i)+a_i-1}{a_i}$. Thus, \eqref{eq:PMnq} transforms into 
\begin{equation}\label{eq:ForMnq}
\mathbb{P}_{\mathcal{M}}( \bfa)  
=  \prod_{k=1}^{n} \frac{1}{k^{a_k} a_k!} \cdot \left( \prod_{i=1}^{n} \prod_{j=0}^{a_i-1} \frac{i\pi_q(i)+ij}{q^i} \right).
\end{equation} 

To connect these probability measures with the derangement probabilities discussed in the introduction, we note that if $\Omega_0$ is the event that $a_1=0$, then 
\begin{equation}\label{eq:Omega0}
P_n = \mathbb{P}_{S_n}(\Omega_0) \quad \mbox{and}\quad P_{n,q} = \mathbb{P}_{\mathcal{M}}(\Omega_0).
\end{equation}

\section{Comparison of probability measures.}\label{proofest}
We prove that the two measures defined above are close in the $\ell^1$-norm.
\begin{theorem}\label{mainlem}
Let $n$ be a positive integer and $q$ a prime power. Then 
\begin{equation}\label{res}
\|\mathbb{P}_{S_n}- \mathbb{P}_{\mathcal{M}}\|_{1} := \sum_{\bfa\in \Omega} |\PP_{S_n}(\bfa)- \PP_{\calM}(\bfa)| = O\Big( \frac{1}{q} \Big).
\end{equation}
\end{theorem}
The bound \eqref{res} immediately implies \eqref{eq:MB} using \eqref{eq:Omega0}.   Indeed, by the triangle inequality, 
\[
|P_n-P_{n,q}| \leq \sum_{\bfa\in \Omega_0}|\PP_{S_n}(\bfa)- \PP_{\calM}(\bfa)|\leq \|\mathbb{P}_{S_n}- \mathbb{P}_{\mathcal{M}}\|_{1} = O\Big(\frac{1}{q}\Big).
\]

The key tool in proving the theorem above is \emph{Gauss's exact formula} for the number of prime polynomials, which may be regarded as comparison of the measures on the event $\{\bfa = (0,\ldots, 0,1)\}$. It is given in terms of  the M\"obius function defined as
\[
\mu(n) = 
\begin{cases}
(-1)^r, & n= p_1\cdots p_r, \mbox{ for distinct prime numbers $p_i$}\\
0, & \mbox{otherwise.}
\end{cases}
\]
Then Gauss's formula is
\begin{equation}\label{Gauss}
i\pi_q(i) = q^i + \sum_{1\neq d\mid i} \mu(d) q^{i/d}.
\end{equation}
The proof is elementary and easy, but beautiful, see \cite[Thm.~2.2]{Rosen}. 
From \eqref{Gauss} one readily derives the useful bounds
\begin{equation}\label{piest}
q^i \ge i\pi_q(i) \ge  q^i - 2\cdot q^{\lfloor i/2 \rfloor}.
\end{equation}

\begin{proof}[Proof of the Theorem]
We may assume that $n>1$ and we put 
\begin{equation}\label{allsum}
X = \|\mathbb{P}_{\mathcal{M}}- \mathbb{P}_{S_n}\|_{1} .
\end{equation}
We write $\mathbb{P}_{\mathcal{M}}(\bfa)$ as
\begin{equation}\label{pnq12}
\mathbb{P}_{\mathcal{M}}(\bfa) = p_{\bfa,1}+ p_{\bfa,2},
\end{equation}
where $p_{\bfa,1} $ and $p_{\bfa,2}$ are the respective contributions from squarefree  and non-squarefree polynomials. Apply the triangle inequality to \eqref{allsum} to obtain  
\[
X \leq \PP(f \mbox{ is not squarefree}) + \sum_{\bfa \in \Omega(n)} | \mathbb{P}_{S_n}(\bfa)- p_{\bfa,1}|.
\]
Here $f$ is sampled uniformly from $\calM$. It is well known that $\PP(f \mbox{ is not squarefree})=\frac{1}{q}$ (see, e.g., \cite[Prop.~2.3]{Rosen}); thus it remains to show
\[
Y = \sum_{\bfa \in \Omega(n)} | \mathbb{P}_{S_n}(\bfa)- p_{\bfa,1}| = O\Big(\frac{1}{q}\Big).
\]
We write $Y=Y_1+Y_2$, according to whether there exists  $j$ with $a_j>\pi_{q}(j)$ or not and show that each $Y_i$ is bounded by $O(1/q)$. To bound $Y_1$, we recall that $a_j$ corresponds to the number of irreducible factors, so the pigeonhole principle tells us the corresponding  polynomial has a repeated factor, hence does not contribute to $p_{\bfa,1}$, so $p_{\bfa,1}=0$. It then follows  by \eqref{eq:ForSn} that
\begin{equation}\label{Y_1bd}
Y_1 \leq 
	\sum_{j=1}^{n}  \sum_{\substack{\bfa\in \Omega(n)\\ a_j>\pi_q(j)}} \mathbb{P}_{S_n}(\bfa) 
	= \sum_{j=1}^{n} \sum_{a_j>\pi_q(j)}  \frac{1}{j^{a_j} a_j!} \sum_{\substack{\bfb \in \Omega(n)\\ b_j=a_j}} \prod_{i \neq j} \frac{1}{i^{b_i} b_i!}
	\leq \sum_{j=1}^{n} \sum_{a_j > \pi_q(j)} \frac{1}{j^{a_j}a_j!},
\end{equation}
where $\sum_{\substack{\bfb \in \Omega(n)\\ b_j=a_j}} \prod_{i \neq j} \frac{1}{i^{b_i} b_i!}\leq 1$ as the probability of a permutation on $n-ja_j$ letters not to have an orbit of size $j$.
Since $\frac{1}{j^{a_{j}+1}(a_j+1)!} \Big/ \frac{1}{j^{a_{j}}a_j!}= \frac{1}{j(a_j+1)}\leq \frac{1}{2}$, we may bound the inner  sum in the right-hand side of \eqref{Y_1bd} by twice the first summand, so by using the lower bound in \eqref{piest} we obtain 
\[
Y_1 \leq 2\sum_{j \ge 1} \frac{1}{j^{\pi_q(j)+1}(\pi_q(j)+1)!} \le 2\sum_{j \ge 1} \frac{1}{(\pi_q(j)+1)!}  = O\Big( \frac{1}{q} \Big).
\]

Now we bound $Y_2$; i.e., considering only $\bfa$'s with $a_j\leq \pi_q(j)$ for all $j$. Similar to the derivation of \eqref{eq:ForMnq}; namely, using $\sum ia_i=n$, we have
\begin{equation}\label{pnq1eval}
p_{\bfa,1} 
	= \frac{1}{q^n}\prod_{i=1}^{n} \binom{\pi_q(i)}{a_i}
	= \prod_{k=1}^{n} \frac{1}{k^{a_k} a_k!}  \prod_{i=1}^{n} \prod_{j=0}^{a_i-1} \frac{i\pi_q(i)-ij}{q^i},
\end{equation}
so 
\begin{equation}\label{defs2}
Y_2=\sum_{\substack{\bfa \in \Omega(n)\\ \forall r : a_r \le \pi_q(r)}} | \mathbb{P}_{S_n}(\bfa)- p_{\bfa,1}|  = \sum_{\substack{\bfa \in \Omega(n)\\ \forall r:a_r \le \pi_q(r)}} \prod_{k=1}^{n} \frac{1}{k^{a_k} a_k!} \left|  1 - \prod_{i=1}^{n} \prod_{j=0}^{a_i-1} \frac{i\pi_q(i)-ij}{q^i} \right|.
\end{equation}
By the upper bound in \eqref{piest} and the assumption $j<a_i\leq \pi_q(i)$, we have $0\leq \frac{i\pi_q(i)-ij}{q^i}\leq 1$. So we may use the Bernoulli-type inequality
\begin{equation}
0 \le 1-\prod_{k=1}^{m}(1-x_k) \le \sum_{k=1}^{m} x_k
\end{equation}
 with $\{ x_k\}_k = \{ 1-\frac{i\pi_q(i)-ij}{q^i} \}_{i,j}$ to obtain
\begin{equation}\label{s2ineq1}
Y_2 
	\le \sum_{\substack{\bfa \in \Omega(n)\\ \forall j : a_j \le \pi_q(j)}} \prod_{k=1}^{n} \frac{1}{k^{a_k} a_k!} \sum_{i=1}^{n} \sum_{j=0}^{a_i-1} ( 1-\frac{i\pi_q(i)-ij}{q^i}).
\end{equation}
From the lower bound in \eqref{piest} we conclude that
\begin{equation}\label{innersums2}
\sum_{j=0}^{a_i-1} ( 1-\frac{i\pi_q(i)-ij}{q^i}) \le 2a_i q^{-\lceil \frac{i}{2} \rceil} + \frac{i}{q^i} \frac{a_i(a_i-1)}{2} \le 4a_i^2 i q^{-\lceil \frac{i}{2} \rceil}.
\end{equation}
Plugging \eqref{innersums2} in \eqref{s2ineq1} yields
\[
\begin{split}
Y_2 
	&\le 4\sum_{\bfa \in \Omega(n)} \prod_{k=1}^{n} \frac{1}{k^{a_k} a_k!} \sum_{i=1}^{n}  a_i^2 i q^{-\lceil \frac{i}{2} \rceil}\\
	& = 4\sum_{i=1}^{n} \sum_{0\leq a\leq n} \sum_{\substack{\bfb \in \Omega(n)\\ b_i=a}} \prod_{\substack{ k=1\\  k \neq i}}^{n} \frac{1}{k^{b_k} b_k!} \frac{1}{i^{a}a!}a^2 i  q^{-\lceil \frac{i}{2} \rceil}\\
	& \le 4\sum_{i=1}^{n} \sum_{0 \le a \le n} \frac{a^2 i q^{-\lceil \frac{i}{2} \rceil}}{i^{a}a!} ,
\end{split}
\]
where as before $\sum_{\substack{\bfb \in \Omega(n)\\ b_i=a}} \prod_{\substack{ k=1\\  k \neq i}}^{n} \frac{1}{k^{b_k} b_k!}  \leq 1$.
Thus,
\[	
Y_2  \le 4\sum_{i \ge 1} i q^{-\lceil \frac{i}{2} \rceil} \sum_{a \ge 0} \frac{a^2 }{a!} =  O\Big( \frac{1}{q}\Big),
\]
as needed to finish the proof.
\end{proof}

\section*{Acknowledgments.}
The authors are partially supported by a grant of the Israel Science Foundation.

%\begin{affil}
%School of Mathematical Sciences, Tel Aviv University, Ramat Aviv, Tel Aviv 6997801, Israel\\
%barylior@post.tau.ac.il
%\end{affil}
%\begin{affil}
%School of Mathematical Sciences, Tel Aviv University, Ramat Aviv, Tel Aviv 6997801, Israel\\
%ofir.goro@gmail.com
%\end{affil}
%\vfill\eject

\end{document}